\documentclass{amsart}
\usepackage{amssymb}
\usepackage{amscd,enumerate,amsfonts,calc,amsmath,verbatim,xypic}         
\newcommand{\Ext}{\operatorname{Ext}}
\newcommand{\TExt}{\operatorname{E\widehat{\vphantom{E}x}t}\!{\vphantom  E}}
\newcommand{\TTor}{\operatorname{T\widehat{\vphantom{T}o}r}\!{\vphantom T }}
\newcommand{\Tor}{\operatorname{Tor}}
\newcommand{\Hom}{\operatorname{Hom}}
\newcommand{\Soc}{\operatorname{Socle}}
\newcommand{\HH}{\operatorname{H}}
\newcommand{\Ima}{\operatorname{Im}}

\newcommand{\T}{{\mathcal T}}
\newcommand{\Coker}{\operatorname{Coker}}
\newcommand{\rank}{\operatorname{rank}}

\newcommand{\depth}{\operatorname{depth}}
\newcommand{\codim}{\operatorname{codim}}
\newcommand{\pd}{\operatorname{pd}}
\newcommand{\vpd}{\operatorname{vpd}}
\newcommand{\vid}{\operatorname{vid}}
\newcommand{\id}{\operatorname{id}}

\newcommand{\col}{\colon}

\newcommand{\ges}{\geqslant}
\newcommand{\les}{\leqslant}
\newcommand{\Ker}{\operatorname{Ker}}

\newtheorem*{Theorem}{Theorem}

\newcommand{\m}{{\mathfrak m}}
\newcommand{\Po}{{\text P}}
\newcommand{\Hi}{{\text H}}
\newcommand{\fm}{{\mathfrak m}}
\newcommand{\fn}{{\mathfrak n}}

\newcommand{\ci}{({\bf ci})}
\newcommand{\gor}{({\bf gor})}
\newcommand{\te}{({\bf te})}
\newcommand{\et}{{(\bf et})}
\newcommand{\ee}{({\bf ee})}
\newcommand{\gap}{({\bf gap})}
\newcommand{\ab}{({\bf ab})}

\newcommand{\bd}{\boldsymbol}
\newcommand{\xra}{\xrightarrow}
\newcommand{\bm}{\begin{matrix}}
\newcommand{\dm}{\end{matrix}}
\theoremstyle{remark}

\swapnumbers
\theoremstyle{plain}
\newtheorem{theorem}{Theorem}[section]

\newtheorem*{Theorem AB}{Theorem AB}

\newtheorem{proposition}[theorem]{Proposition}
\newtheorem{lemma}[theorem]{Lemma}
\newtheorem{corollary}[theorem]{Corollary}
\theoremstyle{definition}
\newtheorem{chunk}[theorem]{}

\newtheorem*{chunk*}{}

\theoremstyle{remark}

\newtheorem{remark}[theorem]{Remark}

\numberwithin{equation}{theorem}
\numberwithin{subchunk}{theorem}
\begin{document}\title[Asymmetric complete resolutions]
{Asymmetric complete resolutions and \\vanishing of Ext over 
Gorenstein rings} 
\author[D.~A.~Jorgensen]{David A.~Jorgensen } 
\address{Department of Mathematics, University of Texas at Arlington, 
Arlington, TX 76019}
\email{djorgens@math.uta.edu}
\author[L.~M.~\c Sega]{Liana M.~\c Sega}
\address{Department of Mathematics, Michigan State University, 
East Lansing, MI 48824}
\curraddr{Department of Mathematics and Statistics, University of
Missouri, Kansas City, MO 64110-2499}
\email{segal@umkc.edu}

\date{\today} 

\begin{abstract}
We construct a class of Gorenstein local rings $R$ which admit
minimal complete $R$-free resolutions $\bd C$ such that  
the sequence $\{\rank_R C_i\}$ is constant for $i< 0$, 
and grows exponentially for all $i>0$.  Over these rings
we show that there exist finitely generated $R$-modules $M$ and $N$ 
such that $\Ext^i_R(M,N)=0$ for all $i> 0$, but $\Ext^i_R(N,M)\ne 0$ 
for all $i>0$. 
\end{abstract}

\maketitle

\section*{Introduction}

Let $R$ be a commutative local Noetherian ring with maximal ideal $\fm$ and 
residue field $k=R/\fm$, and let $M$, $N$ denote finitely generated 
$R$-modules.  We write $\nu_R(M)$ for the minimal number of generators of 
$M$.

It is well-known that $R$ is Gorenstein if and only if the following
remarkable symmetry is satisfied: 
for any module $M$ we have $\Ext^i_R(M,k)=0$ for all $i\gg 0$ if
and only if $\Ext^i_R(k,M)=0$ for all $i\gg 0$ (equivalently,
$\pd_R(M)<\infty$ if and only if $\id_R(M)<\infty$).  A natural
question is whether this statement still holds when $k$ is replaced with any 
module $N$.  More generally, does the Gorenstein
property of $R$ translate into similarities in the asymptotic 
behavior of the sequences 
$\{\nu_R\big(\Ext^i_R(M,N)\big)\}_{i\ge 0}$ and
$\{\nu_R\big(\Ext^i_R(N,M)\big)\}_{i\ge 0}$?

Since complete intersection rings are Gorenstein, 
a foundation is laid by the following theorem of 
Avramov and Buchweitz \cite[5.6]{AB}:

\begin{Theorem AB}
Suppose $R$ is a complete intersection ring.  Then for any pair of 
finitely generated $R$-modules $M$ and $N$ the sequences  
$\{\nu_R\big(\Ext^i_R(M,N)\big)\}_{i\ge 0}$ and 
$\{\nu_R\big(\Ext^i_R(N,M)\big)\}_{i\ge 0}$ both have polynomial growth 
of the {\rm same} degree. 
\end{Theorem AB}

A sequence of positive integers $\{c_i\}_{i\ge 0}$ is said to 
have {\it polynomial growth of degree $d$} if there exist polynomials 
$f(t)$ and $g(t)$, both of degree $d$ and having the same leading term,
such that $g(i)\le c_i\le f(i)$ 
for all $i\gg 0$.  (We adopt the convention that the zero polynomial
has degree $-1$.)  We say that $\{c_i\}_{i\ge 0}$ has 
{\it exponential growth} if there exist $a,b\in\mathbb R$ with  
$1<a\le b$ such that $a^i\le c_i\le b^i$ for all $i\gg 0$.

The numbers
$\beta_i^R(M)=\nu_R\big(\Ext^i_R(M,k)\big)$, with $i\ge 0$, 
are called the {\it Betti\/} numbers of $M$ and the numbers
$\mu^i_R(M)=\nu_R\big(\Ext^i_R(k,M)\big)$, with $i\ge 0$,
are called the {\it Bass\/}
numbers of $M$ (associated to the maximal ideal $\m$). 

Complete intersection rings are characterized by polynomial growth
of the sequence $\{\beta_i^R(k)\}_{i\ge 0}$, see \cite{G1}, \cite{G2},
so Theorem AB does not extend to all Gorenstein rings.  On the other hand, 
it is proved in \cite[8.2.2]{Av2}  
that for every non-complete intersection $R$
the sequence $\{\beta_i^R(k)\}_{i\ge 0}$ has exponential growth.  

We extend Theorem AB to Gorenstein rings with $\m^3=0$ or with
$\codim R\le 4$, where $\codim R=\nu_R(\m)-\dim R$, in the following form:
for every module $M$, both sequences $\{\beta_i^R(M)\}_{i\ge 0}$ and $\{\mu^i_R(M)\}_{i\ge 0}$
have either polynomial growth of the same degree
or exponential growth. However, the
main result of the paper shows that Gorenstein rings are not
characterized by such  symmetry in homological behavior: 

\begin{Theorem}
 There exist Gorenstein rings $R$ with $\fm^4=0$ and $\codim R=6$, and
 finitely generated $R$-modules, such that their Betti
      sequence is constant (respectively, has exponential
      growth)  and their Bass sequence  has
      exponential growth (respectively, is constant). Moreover, there exist finitely generated $R$-modules $M$, $N$
      such that 
$$
\Ext^i_R(M,N)=0\quad\text{for all $i>0$}\quad\text{and}\quad
\Ext^i_R(N,M)\ne 0 \quad\text{for all $i>0$}\,.
$$
\end{Theorem}

The preceeding result answers several questions from the recent 
literature, which we discuss next.

\subsection*{Symmetry in the vanishing of Ext}
Theorem AB, restricted to the case of polynomial growth of degree
$-1$, shows that any complete intersection ring $R$ satisfies the 
following property:
\begin{align*}
{\bf (ee)\quad } &\text{If $M$ and $N$ are finitely generated 
$R$-modules such that $\Ext_R^i(M,N)=0$ }\\
 &\text {for all $i\gg 0$, then $\Ext_R^i(N,M)=0$ for all $i\gg 0$}.
\end{align*}
The authors of \cite{AB} asked whether all Gorenstein 
rings satisfy ${\bf (ee)}$. It was subsequently established in \cite{HJ} 
and \cite{S} that {\bf(ee)} holds for certain classes of Gorenstein local 
rings $(R,\fm)$ other than the complete intersection rings, for example 
Gorenstein rings with $\fm^3=0$, and  Gorenstein rings with 
$\codim R\le 4$.

In \cite{HJ}, Huneke and Jorgensen introduce a class of Gorenstein rings, 
called  AB rings,  and prove that any AB ring satisfies ({\bf ee}). 
In  \cite{JS1}, the present authors constructed Gorenstein rings 
which are not AB, 
but these examples failed to disprove {\bf (ee)}. 

\subsection*{Betti numbers versus Bass numbers; complete resolutions}
The question of whether the Betti and Bass sequences of an $R$-module $M$
have the same asymptotic behavior has been previously posed in the
literature in the more general context of complete resolutions. A {\it
complete resolution} of the $R$-module $M$ is a complex $\bd C$ of
finitely generated free $R$-modules with differentials $d_i\colon
C_i\to C_{i-1}$ such that the complexes $\bd C$ and $\Hom_R(\bd C,R)$
are both exact, and such that $\bd C_{\gg 0}=\bd F_{\gg 0}$ for some
free resolution $\bd F$ of $M$. The complex $\bd C$ is said to be {\it
minimal} if $d_i(C_i)\subseteq \m C_{i-1}$ for all $i\in \mathbb
Z$. If $R$ is Gorenstein, then every finitely generated $R$-module $M$
has a minimal complete resolution $\bd C$. Moreover, any two minimal
complete resolutions of $M$ are isomorphic, cf. \cite[8.4]{AM}, hence
the numbers $\rank_RC_i$ are uniquely determined. We say that $\bd C$
has {\it symmetric growth} if both sequences $\{\rank_R C_i\}_{i\ge
0}$ and $\{\rank_R C_{-i}\}_{i\ge 0}$ have exponential growth or 
polynomial growth of the same degree.

If $R$ is Gorenstein and $\bd C$ is a minimal complete resolution of a
maximal Cohen-Macaulay $R$-module $M$, then $\beta_i^R(M)=\rank_R C_i$
and $\mu_R^{i+d-1}(M)=\rank_R C_{-i}$ for all $i\ge 1$, where $d=\dim
R$. Therefore the question about the Betti and Bass 
sequences can be translated into the question of whether
$\bd C$ has symmetric growth.  Similar versions of this question
have been previously posed in
\cite[9.2]{AM} and in \cite{JS2}.  In \cite{JS2} we constructed 
doubly infinite minimal exact complexes of free modules which had
asymmetric growth; however, the ring was not Gorenstein
and these complexes were not complete resolutions.

\medskip

The paper is organized as follows: In Section 1 we use  results of Avramov \cite{Av} 
and Sun \cite{Su} to prove that
any complete
resolution has symmetric growth over Gorenstein rings with $\fm^3=0$
or $\codim R\le 4$.

In Section 2 we prove that there exist Gorenstein rings $R$ with
$\m^4=0$ and $\codim R = 6$ which admit complete resolutions $\bd C$
for which $\{\rank_R C_i\}_{i\ge 0}$ is constant  
(respectively, grows exponentially), and $\{\rank_R C_{-i}\}_{i\ge 0}$
grows exponentially (respectively, is constant).
We do not know whether such asymmetric complete resolutions  
exist when $\codim R=5$.

Using the results of Section 2, we prove in Section 3  that
there exist finitely generated modules $M$, $N$ which give
counterexamples to ${\bf (ee)}$; the ring is the same as in Section 2.  The 
module $N$ has minimal possible length for such a counterexample,
namely length $2$. 
The results of this section are stated in terms of Tate
(co)homology: we  show that the Tate cohomology groups $\TExt^i_R(M,N)$
vanish for all $i>0$, but do not vanish for all 
$i<0$. 

In Appendix A we establish the structure and relevant 
properties of the rings $R$ from Sections 2 and 3.
These rings are similar to those constructed in \cite {GP}, \cite{JS1}, 
\cite{JS2}.

\section{Symmetric growth of complete resolutions}

In this section we show that
there exist certain classes of Gorenstein rings, other than the class
of complete intersection rings, for which all complete resolutions
have symmetric growth.

Let $(R,\fm,k)$ be a local ring as in the introduction. If $R$ is Gorenstein, then a complete resolution of a finitely
generated $R$-module $M$ has symmetric growth if and only if the Betti
sequence $\{\beta_i^R(M)\}_{i\ge 0}$ and the Bass
sequence $\{\mu^i_R(M)\}_{i\ge 0}$ have the same growth. 
For the convenience of the reader, we prove this in Lemma \ref{CRandBB}.

\begin{chunk}
\label{MCM}
The asymptotic behavior of the Betti sequences remains
unchanged upon passing to syzygies. When the ring is Gorenstein, the
same is true for the Bass sequences. We 
may  thus assume, whenever convenient, that $M$ is a maximal 
Cohen-Macaulay module over $R$. 
\end{chunk}

We let $M^*$ denote the $R$-module $\Hom_R(M,R)$. If $\bd D$ is a complex,
then $\bd D^*$ denotes the complex with $(\bd D^*)_i=(D_{-i})^*$ for each
$i$, and with induced  differentials.

\begin{lemma}\label{CRandBB} Let $R$ be a Gorenstein local ring
of dimension $d$, let $M$ be a finitely generated maximal Cohen-Macaulay 
$R$-module, and $\bd C$ a minimal complete resolution of $M$.  The
following equalities hold:
\begin{enumerate}[\quad\rm(1)]
\item $\beta_i^R(M)=\rank_R C_i$ for all $i\ge 0$, and
$\mu_R^{i+d-1}(M)=\rank_R C_{-i}$ for all $i\ge 1$;
\item $\beta_i^R(M^*)=\mu_R^{i+d}(M)$ for all $i\ge 0$.
\end{enumerate}
\end{lemma}

\begin{proof} (1) 
If $d=0$, the statement is clear: $\bd C_{\ges 0}$ is 
a minimal free resolution of $M$ over $R$, and $\bd C_{\les -1}$ is a 
minimal injective resolution of $M$ over $R$. 
If $d>0$, then let  $\bd
x=x_1,\dots,x_d$ be a maximal regular sequence for both $R$ and
$M$. Note that 
$\bd C/(\bd x)\bd C=\bd C\otimes_R R/(\bd x)$ 
is a minimal complete resolution of $\overline M=M/(\bd x)M$ 
over the zero-dimensional Gorenstein ring $\overline R=R/(\bd x)$. The
conclusion then follows from the isomorphisms 
$\Ext_R^i(M,k)\cong\Ext_{\overline R}^i(\overline M,k)$ 
and 
$
\Ext_R^{i+d}(k,M)\cong\Ext_{\overline R}^i(k,\overline M)
$, which hold for for all $i\ge 0$, cf. for example \cite[p. 140]{M}.

(2) Note that  $(\bd C_{\les -1})^*$ is a minimal free resolution of $M^*$ over
$R$, hence for all $i\geq 0$ we have $\beta_i^R(M^*)=\rank_R (C_{-i-1})^*=
\rank_R C_{-i-1}$. By part (1), the last expression is equal  $\mu_R^{i+d}(M)$.
\end{proof}

\begin{chunk}
\label{minimal multiplicity} 
If $R$ is a Gorenstein ring with $\codim R\ge 2$, then
$R$ has multiplicity at least $\codim R+2$. Otherwise, when $\codim R\leq 1$,
the multiplicity is at least $\codim R+1$.  In either case, when
equality holds we say that $R$ is {\em Gorenstein of minimal 
multiplicity}. If  $R$ is furthermore Artinian, then $R$ has minimal
multiplicity if and only if $\fm^3=0$. 
\end{chunk}

\begin{chunk}
\label{multiplicity} If $R$ is Gorenstein of minimal multiplicity
with $\codim R\ge 3$, then for each finitely generated $R$-module $M$ 
either $M$ has finite projective dimension, or the
sequence $\{\beta^R_i(M)\}_i$ has exponential growth.
Indeed, if $\dim R=0$, then $\fm^3=0$ and the result is proved by 
Sj\"odin \cite{Sj}, cf. also Lescot \cite{L}.  If $\dim R>0$, then 
the reduction to the zero dimensional case can be done  as described in \cite[1.7]{S}.
\end{chunk}

\begin{chunk}
\label{codim}
Assume now that $R$ is Gorenstein and $\codim R\le 4$. Avramov \cite{Av}
and Sun \cite{Su} classified the possible behavior of the Betti numbers 
of a finitely generated $R$-module $M$. They show that the Betti sequence 
has either polynomial growth, or exponential growth. Note that Avramov 
and Sun use the terminology  of {\it strong} polynomial/exponential growth 
for describing the same concepts that we are concerned with, only that we 
omit the word ``strong''.  The classification involves the notion
of {\it virtual projective dimension\/}.  We recall this notion for the 
reader's convenience:  
let $M$ be a finitely generated module over a local ring $R$ 
(not necessarily Gorenstein).  If the residue field $k$ of $R$ is 
infinite, set $\widetilde R=\widehat R$, the $\m$-adic completion of $R$;
if $k$ is finite, set $\widetilde R$ to be the maximal-ideal-adic
completion of $R[Y]_{\m R[Y]}$, where $Y$ is an indeterminate.
We say that a map of local rings $\widetilde R\gets Q$ is an {\it
embedded deformation\/} of $R$ if its kernel is generated by a $Q$-regular
sequence contained in the square of the maximal ideal of $Q$.   
The {\it virtual projective dimension\/} of $M$ is the number
$$
\vpd_R M=\min\{\pd_Q (M\otimes_R \widetilde R) \,\vert\, \widetilde R\gets
Q \text{ is an embedded deformation of $\widetilde R$}\}
$$
A similar invariant, called {\it virtual injective dimension} and denoted
$\vid_RM$, can be defined by replacing $\pd_Q (M\otimes_R\widetilde R)$
with $\id_Q (M\otimes_R\widetilde R)$ in the formula above. 
In \cite{Av} and \cite{Su} it is shown that if $R$ is Gorenstein
with $\codim R\le 4$, then the Betti numbers
of $M$ fall into one of two categories, each described by equivalent
conditions:
\begin{enumerate}
\item $\{\beta_i^R(M)\}_i$ has polynomial growth if and only if
$\vpd_R M<\infty$; when $\vpd_R M$ is finite it is equal to
$\vpd_R M = \depth R-\depth M+q+1$, where $q$ is the degree of 
polynomial growth of the sequence;
\item $\{\beta_i^R(M)\}_i$ has exponential growth if and only if
$\vpd_R M = \infty$.  
\end{enumerate}
Using Lemma \ref{CRandBB}(2), we see that the Bass numbers of $M$
have the same behavior: they either have polynomial growth or 
they have exponential growth. 
\end{chunk}

We are now ready to prove the main result of this section, which
assembles the results of \cite{Av} and \cite{Su} listed above.

\begin{theorem}
\label{BandB}
Assume that $R$ is Gorenstein, and either $R$ has minimal 
multiplicity, or $\codim R\le 4$. If $M$ is a finitely generated
$R$-module, then one of the following statements is satisfied:
\begin{enumerate}[\quad\rm(1)]
\item Both sequences $\{\beta_i^R(M)\}_{i\ge 0}$ and
      $\{\mu^i_R(M)\}_{i\ge 0}$ have  polynomial growth of
      the same degree. 
\item Both sequences $\{\beta_i^R(M)\}_{i\ge 0}$ and
      $\{\mu^i_R(M)\}_{i\ge 0}$ have  exponential growth.
\end{enumerate}
\end{theorem}

\begin{proof}
Assume first that $R$ has minimal multiplicity.
If $\codim R\le 2$, then $R$ is a complete intersection, 
and Theorem AB in the introduction  shows that $M$ satisfies (1). 
Assume now $\codim(R)\ge 3$. If $\pd_R M=\infty$, then $\id_RM=\infty$ as
well, and  it can be 
immediately seen from \ref{CRandBB}(2) 
and \ref{multiplicity} that $M$ satisfies condition (2). If $\pd_RM$
is finite, then $\id_RM$ is also finite, and $M$ satisfies (1). 

Assume $\codim R\le 4$.  Since
 $R$ is Gorenstein, the same is true for any ring $Q$ in  an embedded deformation 
$\widetilde R \gets Q$, and so it is clear that 
 $\vid_RM=\vpd_RM+\depth M$. It follows then
 from \cite[1.8]{Av} and \ref{codim} that $M$ satisfies (1)
 whenever $\vpd_RM$ and $\vid_RM$ are finite and $M$ satisfies (2)
 whenever they are infinite. 
\end{proof}

Lemma \ref{CRandBB} and \ref{MCM} now yield:

\begin{corollary} Assume $R$ is Gorenstein. If  $R$ has minimal 
multiplicity, or if $\codim R\le 4$, then any minimal complete
resolution over $R$ has symmetric growth. \qed
\end{corollary}

\section{Asymmetric growth of complete resolutions}

In this section we show that complete resolutions need not be
symmetric when $\m^4=0$ and $\codim R=6$.

Let $k$ be a field which is not algebraic over a finite field and let
$\alpha\in k$ be an element of infinite multiplicative
order. Throughout the whole section we consider the ring $R$ to be
defined as follows. 

\medskip

\begin{chunk}
\label{define}
Let $P=k[T,U,V,X,Y,Z]$ be the polynomial ring in six variables 
(each of degree one) and set 
$R=P/I$,
  where $I$ is the ideal generated by the following fifteen quadratic
  polynomials:
\begin{gather}
Z^2,\,\,UZ-TX-\alpha UV,\,\,U^2,\,\,YZ+VY,\,\,UY,\,\,
Y^2-TX-(\alpha-1)UV,\notag\\
XZ+\alpha VX,\,\,UX,\,\,XY,\,\,X^2-TX-TV,\,\,
TZ+TY+\alpha VX,\,\,TU,\notag\\
TY-VX+TV,\,\,T^2+(\alpha +1)UV-VY,\,\,V^2.
\notag
\end{gather}
\end{chunk}
  
Let $t,u,v,x,y,z$ denote the residue classes of the 
variables modulo $I$, and $\fm$ denote the ideal they 
generate.

\begin{proposition} 
\label{prop} The ring $R$ is local, with maximal ideal $\fm$, and satisfies the following properties: 
\begin{enumerate}[\quad\rm(1)]
\item $R$ has Hilbert series $H_R(t)=1+6t+6t^2+t^3$. More precisely, a basis of $R$ over $k$ is given by the following 
fourteen elements:
       $$
1,\, t, \, u,\, v,\, x,\, y,\, z,\, tv,\,uv,\,vx,\,vy,\,vz,\,tx,\,vtx\,
$$
\item $R$ is Gorenstein, with $\Soc(R)=(tvx)$. 
\item $R$ is a Koszul algebra. 
\end{enumerate}
\end{proposition}

We prove this proposition in the Appendix as Theorem A.1; 
we also provide there a relevant part of the multiplication table of $R$. 

\begin{chunk} For each $i\leq 0$ we let $d_i\colon R^2\to R^2$ 
denote the map given 
with respect to the standard basis of $R^2$ by the matrix
\[
\left(\bm
v& y\\
\alpha^{1-i}x& z\dm\right).
\]
Let $d_1\col R^3\to R^2$ denote the map represented with 
respect to the standard bases of $R^3$ and $R^2$ by the matrix 
$$
\left(\begin{matrix}v&y&0\\
x&z&tv\end{matrix}\right).
$$ 
Consider a minimal free resolution of $\Coker d_1$ with $d_1$ as
the first differential:
$$
\cdots \to R^{\beta_i}\xrightarrow{d_i} R^{\beta_{i-1}}\to \cdots\to
R^{\beta_2}\xrightarrow{d_2}R^3\xrightarrow{d_1}R^2
$$ 
\end{chunk}

\begin{theorem}
\label{crazy}
The sequence of homomorphisms 
$$ 
\bd C:\quad\cdots\to R^{\beta_2}\xra{d_2}
R^3\xra{d_1}R^2\xra{d_0}R^2\xra{d_{-1}}R^2
\xra{d_{-2}} R^2\to\cdots 
$$ 
is a minimal complete resolution such the following hold:
\begin{enumerate}[\quad\rm(1)]
\item The sequence $\{\rank C_i\}_{i\geq 0}$ has exponential growth.
\item $\rank C_i=2$ for all $i\le 0$. 
\end{enumerate}
\end{theorem}
\begin{proof}[Proof of (2).] We postpone the proof of (1) 
to the end of the section. 
The minimality of $\bd C$
is clear from the definition of the differentials $d_i$.  
Moreover, the defining 
equations of $R$ guarantee $d_id_{i+1}=0$ for all $i\leq 0$, 
hence $\bd C$ is a complex.  Since the ring $R$ is
Gorenstein, the exactness of $\bd C$ also implies that the complex 
$\Hom_R(\bd C, R)$ is exact.  Therefore it remains to show  that 
$\bd C$ is exact.   

 Let $(a,b)$ denote an element of $R^2$ written in the standard
basis of $R^2$ as a free $R$-module. One may check that for each 
$i\le 0$ the $k$-vector
space $\Ima d_{i}$ has the following fourteen 
linearly independent elements:
\begin{align*}
d_i(1,0)&=(v,\alpha^{1-i}x)      &d_i(tv,0)&=(0,\alpha^{1-i}tvx)\\
d_i(t,0)&=(tv,\alpha^{1-i}tx)    &d_i(tx,0)&=(tvx,0)\\  
d_i(u,0)&=(uv,0)                 &d_i(0,1)&=(y,z)\\
d_i(v,0)&=(0,\alpha^{1-i}vx)     &d_i(0,t)&=(vx-tv,tv-(\alpha+1)vx)\\
d_i(x,0)&=(vx,\alpha^{1-i}(tv+tx))    &d_i(0,u)&=(0,\alpha uv+tx)\\
d_i(y,0)&=(vy,0)                 &d_i(0,v)&=(vy,vz)\\
d_i(z,0)&=(vz,-\alpha^{2-i}vx)   &d_i(0,y)&=((\alpha-1)uv+tx,-vy).
\end{align*}

We thus have $\rank_k \Ima d_i\ge 14$ for all $i\le 0$. Since $\rank_k
R^2=28$, it follows that $\rank_k \Ker d_i\le 14$ for all $i\le 0$. 
We conclude that $H_i(\bd C)=0$ for $i\leq -1$ and 
$\rank_k\Ker d_0=14$ 
 
For $i=1$ the images above are also those of $d_1$, by replacing 
$d_i(1,0)$ with $d_1(1,0,0)$ and so on. Here $(a,b,c)$  denotes an 
element of $R^3$ in its standard basis as a free $R$-module. 
Not more than thirteen of these images are  linearly independent, since we 
have the relation
$$
(vx,tv+tx)=(tv,tx)+(vx-tv,tv-(\alpha+1)vx)+(\alpha+1)(0,vx)\,.
$$
One can check that we are left with precisely thirteen linearly
independent elements, and a fourteenth element in $\Ima d_1$ can 
be chosen to be $d_1(0,0,1)=(0,tv)$. 
Thus $\rank_k \Ima d_1\ge 14$ and it follows that 
$H_0(\bd C)=0$. The definition of $\bd C_{\ge 1}$ gives that
$H_i(\bd C)=0$ for all $i\ge 1$, and therefore $\bd C$ is
exact. 
\end{proof}

Next we provide the necessary background for the proof of part (1) of 
the Theorem.   

\begin{chunk}
\label{maps}
Let $\pi\col A\to B$ be a  ring homomorphism, $D$ an $A$-module, $E$
a $B$-module (with the $A$-module structure induced by $\pi$) and $\phi\col
D\to E$ a homomorphism of $A$-modules. Then for each $B$-module
$L$ one has a natural homomorphism $$\Tor^{\pi}(L,\phi)\col 
\Tor^A(L,D)\to \Tor^B(L,E)$$ which may be computed as follows:
Let  $\bd D$ be a free $A$-resolution of $D$ and $\bd E$ a free $B$-resolution
of $E$. Let $\bd{\widetilde \phi}\col B\otimes_A\bd D\to\bd E$ be a
lifting of $\phi$ to a homomorphism of complexes of $B$-modules. 
The homomorphism
$\Tor^{\pi}(L,\phi)$ is then induced in homology by the following
homomorphism of complexes, which is unique up to homotopy: 
$$
L \otimes_A\bd D=L\otimes_B(B\otimes_A\bd D)
\xrightarrow{L\otimes_B\bd{\widetilde \phi}}L\otimes_B\bd E
$$
\end{chunk}
We say that a positively graded ring
$A=\bigoplus_{i\ge 0}A_i$ is {\it standard graded} if $A_0$ is a 
field and $A$ is generated over $A_0$ by $A_1$.

\begin{proposition}
\label{ap2}
Let $A$ be a standard graded algebra over a field $\ell$ and set
$\fn=\bigoplus_{i\ge 1}A_i$. Let $\pi\colon A\to B$ be a 
surjective homomorphism of graded
rings with $\Ker \pi\subseteq \fn^2$.  
Assume that $D$ is a finitely generated graded $A$-module with a linear 
graded free resolution 
and that $E$ is a finitely generated graded $B$-module. 

Suppose  
there exists a
 homomorphism of $A$-modules $\phi\colon D\to E$ such that the induced
 map $\overline \phi\col D/\fn D\to E/\fn E$ is injective. Then the induced
 homomorphisms 
$$\Tor^{\pi}_i(\ell,\phi)\col \Tor^A_i(\ell,D)\to \Tor^B_i(\ell,E)$$ 
are injective for each $i$. 
\end{proposition}

\begin{proof}

Consider a linear resolution $(\bd D,\delta)$ of $D$, together with an
augmentation map $\varepsilon\col D_0\to D$.  We will construct
inductively a minimal graded free resolution $\bd E$ of $E$, with an
augmentation map $\eta\col E_0\to E$, and a map of complexes of
$A$-modules $\bd \varphi\col \bd D\to \bd E$ such that the induced map
$B\otimes_A\bd \varphi\col B\otimes_A\bd D\to\bd E$ is a split
injection in each homological degree.

Let $e_1, \dots, e_{a}$ denote the standard basis of $D_0=A^{a}$. 
Then $\varepsilon(e_1),\dots ,\varepsilon (e_{a})$ form a homogeneous 
minimal system of generators for $D$. Since the induced map
$\overline \phi \col D/\fn D\to E/\fn E$ is injective, the elements 
$\phi(\varepsilon(e_1)),\dots,\phi(\varepsilon (e_{a}))$ are part of  
a minimal system of generators
for $E$. This shows that we can choose $E_0=B^b$, with $b\ge a$, and the
map $\eta \col B^{b}\to E$ can be chosen so that
$\eta(f_i)=\phi(\varepsilon(e_i))$ for each $i$ with $1 \le i\le a$,
where  $f_1,\dots,f_{b}$ is  the standard basis of $B^{b}$.

 If we define an $A$-module homomorphism $\varphi_0\col A^{a}\to B^{b}$
such that $\varphi_0(e_i)=f_i$, then the right-hand part of the diagram
below is commutative. We set $D'=\Ker(\varepsilon)$ and
$E'=\Ker(\eta)$, and we let $\phi\,'\col D'\to E'$ denote the induced
$A$-module homomorphism, which makes entire diagram commutative: 
\[
\xymatrixrowsep{1.6pc} \xymatrixcolsep{0.4pc} 
\xymatrix{
0\ar@{->}[rr]& &D'\ar@{->}[rrr]\ar@{->}^{\phi'}[d]
& & & A^{a}\ar@{->}[rrr]^{\varepsilon}\ar@{->}[d]^{\varphi_0}
& & &D\ar@{->}[d]^{\phi}\ar@{->}[rr]& &0\\
0\ar@{->}[rr]& &E'\ar@{->}[rrr]& & &B^{b}
\ar@{->}[rrr]^{\eta}& & &E\ar@{->}[rr]& &0
}
\]
We want to prove that the induced map 
$\overline{\phi'}\col D'/\fn D'\to E'/\fn E'$ is injective. 

Let $g_1,\dots,g_{a'}$ be the standard basis of $D_1=A^{a'}$. 
The elements $\delta_1(g_1),\dots,\delta_1(g_{a'})$ form a minimal 
system of generators for $D'$. Let $\alpha_i\in A$ be such that 
$$
\varphi_0\big(\alpha_1\delta_1(g_1)+\dots+
\alpha_{a'}\delta_1(g_{a'})\big)\in \fn E'\subseteq \fn^2B^{b}.
$$
Since the matrix representing $\delta_1$ has linear entries, 
we can think of $\delta_1(g_i)$ as column vectors with components in
$A_1$ and hence of
$\varphi_0(\delta_1(g_i))$ as column vectors with components in
$B_1$. Thus, the degree one part of 
the above expression is equal to zero, so we get
$$
\varphi_0\left(\sum_i\overline \alpha_i\delta_1(g_i)\right)=
\sum_i\overline \alpha_i\varphi_0(\delta_1(g_i))=0
$$
where $\overline\alpha_i$ denotes the degree zero component 
of $\alpha_i$.  Since the homomorphism $\pi\col A\to B$ has 
$\Ker \pi \subseteq \fn^2$, 
we conclude that $\sum_i\overline \alpha_i\delta_1(g_i)\in \fn^2A^{a}$, 
hence, by degree considerations,
$\sum_i\overline \alpha_i\delta_1(g_i)=0$. Since the elements 
$\delta_1(g_i)$ in this sum are part of a homogeneous minimal system of 
generators for $D'$, it follows that
$\overline\alpha_i=0$ for all $i$. Therefore $\alpha_i\in \fn$ for all 
$i$, and this  shows $\overline{\phi'}$ is injective.  
 
Using the construction above as the induction step, we obtain then a
resolution $\bd E$ and a homomorphism of complexes $\bd \varphi\col\bd
D\to \bd E$.  The homomorphism of complexes $\bd{\widetilde \phi}=\bd
\varphi\otimes_AB\col \bd D\otimes_AB\to \bd E$ is then a lifting of
$\phi$ and is a split injection in each degree. This gives the desired 
conclusion. 
\end{proof}

Recall that $P$ denotes the polynomial ring $k[T,U,V,X,Y,Z]$ and that the
monomials $U^2$ and $UY$ are among the generators of the ideal $I$
defining $R$ as $P/I$.  

\begin{lemma}
\label{A}
Consider the ring $A=P/(U^2,UY)$, the $A$-module $D=(U,Y)A$, and the 
$R$-module $E=\ker d_1$.
Let $\pi \col A\to R$ denote the canonical
projection. The following then hold:
\begin{enumerate}[\quad\rm(1)]
\item The $A$-module $D$ has a linear resolution and its
      Poincar\'e series is equal to $(2+t)(1-t-t^2)^{-1}$.
\item There exists a homomorphism of $A$-modules
$\phi\col D\to E$ such that the induced map 
$\Tor^{\pi}_i(k,\phi)\col
\Tor_i^A(k,D)\to\Tor_i^R(k,E)$
is injective for each $i$. 
\
\end{enumerate}
\end{lemma}

\begin{proof}
(1) Set 
$Q=k[U,Y]/(U^2,UY)$ and let $\bd G$ denote a minimal free resolution
of the residue field $Q/(U,Y)Q$ over $Q$. Note that $A\cong
Q\otimes_kk[T,V,X,Z]$, and a minimal free resolution of $A/D$
over $A$ is given by the complex $\bd D=\bd G\otimes_kk[T,V,X,Z]$. 

Since  $Q$ is a Koszul algebra (see for example \cite{F}),
the resolution $\bd G$ is linear, and the Poincar\'e series of
$A/D$ over $A$ is
\[
\Po^A_{A/D}(t)=\Po^Q_{Q/(U,Y)Q}(t)=\frac{1}{\Hi_Q(-t)}=\frac{1+t}{1-t-t^2}\,.
\]
The Poincar\'e series of $D$ is then equal to
$t^{-1}(\Po^A_{A/D}(t)-1)$. 

(2) Set $p=(0,0,u)$ and $q=(0,0,y)$, considered as elements of $R^3$.  
It can be easily
checked that $p,q\in \Ker d_1$, hence $Rp+Rq\subseteq E$. We define 
$\phi\col D\to E$ as the following composition:
$$
\phi\col D\xra{} Rp+Rq\hookrightarrow E,
$$ 
where the leftmost map is the restriction of the map $\varphi\col A\to
R^3$ given by $\varphi(r)=\big(0,0,\pi(r)\big)$. Note that 
$ap+bq\in \fm E$ for some $a,b\in R$ implies 
$ap+bq\in\fm^2R^3$, and hence $a,b\in \fm$ by degree considerations.
This shows that the induced map $\overline \phi\col D/\fn D\to E/\fn E$ is
injective, where $\fn$ denotes the maximal ideal of $A$. We can then
apply Proposition \ref{ap2}. 
\end{proof}

\begin{proof}[Proof of Theorem {\rm  \ref{crazy}(1)}]
We use the notation in the statement of the Lemma above.
Part (1) of the Lemma shows that the sequence 
$\{\rank_k\Tor_i^A(k,D)\}_i$ has exponential growth. From part (2) 
we conclude that the sequence $\{\rank_k\Tor_i^R(k,E)\}_i$ has 
exponential growth, as well. Note that a minimal free resolution 
of the $R$-module $E$ is given by the truncation $\bd C_{\ges 2}$, 
hence $\rank_k C_{i+2}=\rank_k\Tor_{i}^R(k,E)$ for all $i\geq 0$.
\end{proof}

\section{Asymmetry in the vanishing of Ext}

We will use the notation introduced in the second section. 
In particular, the ring $R$ is the one defined in \ref{define}. 
Recall that $R$ is zero-dimensional and Gorenstein.

\begin{chunk}
\label{tate}
If $\bd T$ is a complete resolution of the $R$-module $X$, 
and $Y$ is an $R$-module,
then for each $i$ the {\it Tate (co)homology} groups are defined by
$$
\TExt^i_R(X,Y)=\HH_{-i}\Hom(\bd T, Y)\quad\text{and}
\quad\TTor_i^R(X,Y)=\HH_i(\bd T\otimes_RY)\,.
$$
Since the ring $R$ is zero-dimensional Gorenstein, the complete
resolution $\bd T$ can be chosen to agree with a minimal free 
resolution of $X$ in all
nonnegative degrees, cf. \cite[3.1]{AM}, for example.  Thus for all 
$i>0$ there are isomorphisms
$$
\TExt^i_R(X,Y)\cong\Ext^i_R(X,Y)\quad\text{and} 
\quad\TTor_i^R(X,Y)\cong \Tor_i^R(X,Y)\,.
$$
Also, for all $i$ one has
$$
\TExt^{-i-1}_R(X,Y)\cong \TTor_i^R(X^*,Y)\,.
$$
Matlis duality yields for all $i$ the following isomorphisms:
$$
\Tor_i^R(Y,X^*)^*\cong \Ext^i_R(Y,X).
$$
\end{chunk}

\medskip

Recall that $d_i$ denotes the differential of the complex $\bd C$ 
defined in Section 2.

\begin{theorem}
\label{eedies}
Set $M=\Coker d_0^*$ and $N=R/(t,u,v-x,y-x,z-x)$. The following then hold:
\begin{enumerate}[\quad\rm(1)]
\item $\TExt_R^i(M,N)=0$ for all $i>0$ 
\item $\TExt_R^i(M,N)\neq 0$ for all $i<0$.
\item $\TTor_i^R(M,N)=0$ for all $i>0$.
\item $\TTor_i^R(M,N)\ne 0$ for all $i<0$. 
\end{enumerate}
\end{theorem}

In view of the isomorphisms in \ref{tate}, we conclude  the Gorenstein
ring $R$ does not have the property
{\bf (ee)}: 

\begin{corollary}
\label{eecor}
For $R$, $M$ and $N$ as above we have:
\begin{align*}
&\text{$\Ext^i_R(M,N)=0$  for all $i>0\,;$}\\
&\text{$\Ext^i_R(N,M)\ne 0$ for all  $i>0$.}
\end{align*}
\end{corollary}

The proof of the Theorem is given after the following lemma, which 
contains one of the ingredients of the proof.

\begin{lemma}
\label{length2}
Set $E=\ker  d_1$. If $L$ is an $R$-module of length two, then the 
following hold:
\begin{enumerate}[\quad\rm(1)]
\item $\Tor^R_i(E,L)\ne 0$ for all $i>0$. 
\item $\Ext_R^i(E,L)\ne 0$ for all $i>0$. 
\end{enumerate}
\end{lemma}

\begin{proof}
(1) Since $L$ has length two, there is an exact sequence
$$
0\to k\to L\to k\to 0.
$$

Using the notation of Lemma \ref{A}, and the naturality of the maps 
defined in \ref{maps}, the short exact sequence above induces long 
exact sequences both over $R$ and over $A$, and they can be embedded  
in a commutative diagram as follows: 
\[
\xymatrixrowsep{1.6pc} \xymatrixcolsep{1.4pc} 
\xymatrix{
\cdots \ar@{->}[r]& \Tor_i^A(k,D)\ar@{->}[r]
\ar@{^{(}->}^{\Tor^{\pi}_i(k,\psi)}[d]&\Tor_i^A(L,D)
\ar@{->}[r]\ar@{->}[d]^{\Tor^{\pi}_i(L,\psi)}&\Tor_i^A(k,D)
\ar@{^{(}->}[d]^{\Tor^{\pi}_i(k,\psi)}
\ar@{->}[r]^{\Delta_i^A}&\Tor_{i-1}^A(k,D)\ar@{->}[r]
\ar@{^{(}->}[d]^{\Tor_{i-1}^{\pi}(k,\psi)}&\cdots\\
\cdots \ar@{->}[r]&\Tor_i^R(k,E)\ar@{->}[r]&\Tor_i^R(L,E)
\ar@{->}[r]&\Tor_i^R(k,E)\ar@{->}[r]^{\Delta_i^R}&\Tor_{i-1}^R(k,E)
\ar@{->}[r]&\cdots
}
\]
Counting from the left we have: the first, third and fourth vertical 
maps are injective, cf. Lemma \ref{A}(2).  
If $\Tor_i^R(E,L)=0$ for some $i>0$, then the connecting 
homomorphism $\Delta_i^R$ is injective, so the commutativity of the 
rightmost square implies that $\Delta^A_i$ is injective. 
However, this is not possible, because Lemma \ref{A}(1) shows that 
the Betti numbers of $D$ over $A$ are strictly increasing. 

(2) By Matlis duality, the $R$-module $L^*$ has length two. We can 
then apply part (1) to conclude $\Tor_i^R(E,L^*)\ne 0$ for all $i>0$ 
and then use  the isomorphism $\Ext^i_R(E,L)^*\cong \Tor_i^R(E,L^*)$.  
\end{proof}

\begin{proof}[Proof of Theorem {\rm \ref{eedies}}]
(1) A complete resolution of $M$ is given by the complex
$\Hom_R(\bd C,R)$. We have 
$$
\TExt^i_R(M,N)=\HH_{-i}(\Hom_R(\Hom_R(\bd C,R),N))\cong \HH_{-i}(\bd
C\otimes_R N)
$$
In negative degrees $\bd C\otimes_R N$ is  the complex
\[
 N^2
\xra{\bigl(\begin{smallmatrix} x&x\\\alpha x&x\end{smallmatrix}\bigr)}N^2
\xra{\bigl(\begin{smallmatrix} x&x\\\alpha^2 x&x\end{smallmatrix}\bigr)}N^2
\xra{\bigl(\begin{smallmatrix} x&x\\\alpha^3 x&x\end{smallmatrix}\bigr)}
N^2\to\cdots
\]
Since $N\cong k[x]/(x^2)$, this complex is acyclic, hence 
$\TExt_R^i(M,N)=0$ for all $i>0$.

(2)  By the isomorphisms in \ref{tate} we need to show that 
$\Tor_i^R(M^*,N)\ne 0$ for
    all $i> 0$. Note that $M^*=\Coker d_2=\Ima d_1$, hence the module $E$  in
    Lemma \ref{length2} is the first syzygy of $M^*$. Since $N$ has 
length $2$, the Lemma shows
    $\Tor_i^R(M^*,N)\ne 0$ for all $i\ge 2$. To show that
    $\Tor_1^R(M^*,N)\ne 0$, consider the short exact sequence
$$
0\to k\to N\to k\to 0.
$$
and the induced long exact sequence
$$
\Tor_1^R(M^*,N)\to\Tor_1^R(M^*,k)\to M^*\otimes_Rk\to M^*\otimes_RN\to
M^*\otimes_Rk\to 0.
$$
If $\Tor_1^R(M^*,N)=0$, then $\beta_2=\dim_k\Tor_1(M^*,k)\le
\dim_k (M^*\otimes_Rk)=3$.  However, this contradicts the fact that
the following five elements are part of a 
minimal system of generators for $\Ker d_1$: 
$$
(0,0,t),\, (0,0,u),\, (0,0,v),\, (0,0,y),\, (0,0,z).
$$

The proofs of (3) and (4) are similar. 
\end{proof}

\begin{remark}
\label{aus}
The following conjecture of Auslander appears in \cite[p.\ 795]{A} and 
\cite{H}: 
{\it let $\Lambda$ be an Artin algebra.  For every finitely generated
$\Lambda$-module $M$ there exists an integer $n_M$ such that for all
finitely generated $\Lambda$-modules $N$, if
$\Ext^i_\Lambda(M,N)=0$ for all $i\gg 0$, then 
$\Ext^i_\Lambda(M,N)=0$ for all $i\ge n_M$.}

Original counterexamples to this conjecture were given
by the authors in \cite{JS1} over a class of codimension-five 
Gorenstein rings. The codimension-six Gorenstein rings $R$ 
of this paper provide other counterexamples, as follows.

For each $q\ge 1$ set  
$$
N_q=R/(t,u,v-\alpha^qx,v-y,v-z).
$$
A computation similar to that in \cite{JS1} shows that the following
holds:
\[
\Ext^i_{R}(M,N_q)=0\quad \text{if and only if \ $i\ne 0,q-1,q$}.
\]

Subsequent to \cite{JS1}, Smal{\o} \cite{Sm}, and independently Mori
\cite{Mo}, gave a very simple non-commutative counterexample to the 
conjecture of Auslander, involving the ring
$A=k\langle x,y\rangle/\langle x^2,y^2,xy-\alpha yx\rangle$, where 
$\alpha\in k$ has infinite multiplicative order. However, this ring fails 
to supply counterexamples to $\ee$.  Indeed, all the 
indecomposable modules over $A$ have been classified, see for example
\cite[Section 4]{Sc}.  They are $A$, syzygies and cosyzygies of $k$, 
modules of periodicity one, and non-periodic modules having bounded 
resolutions.  The differentials in the resolutions of the latter type 
are described as follows: there exists a square matrix $M(t)$ 
with entries in $A[t]$, with $t$ an indeterminate, such that the $i$th 
differential of the resolution is represented by $M(\alpha^i)$. 
An argument similar to that of \cite[3.13]{JS1} 
may be used to see that the ring $A$ satisfies $\ee$. 
\end{remark}

\begin{remark}
The paper \cite{AB} initiated the question of whether
there exist homologically defined classes of
rings strictly contained between the complete intersections and
the Gorenstein rings.  The results of \cite{HJ} and \cite{JS1} show
that the answer is positive: one such class consists of all 
rings $R$ with the property $\ab$ that $R$ is Gorenstein and there exists an integer $n$, depending
only on the ring, such that for all pairs $(M,N)$ of 
finitely generated $R$-modules one has 
$$
\text{$\Ext^i_R(M,N)=0$ for all $i\gg 0$ implies $\Ext^i_R(M,N)=0$ for all
$i> n$.}
$$
The class of Gorenstein rings satisfying $\ab$ was studied in \cite{HJ}; 
it is proved there that this class properly contains the 
complete intersections. The counterexamples to Auslander's 
conjecture in \cite{JS1}, and those in this paper, show that the
Gorenstein rings satisfying $\ab$ are properly contained in the class
of Gorenstein rings. 

Other classes of Gorenstein rings previously considered are those
defined by the properties:
\begin{align*}
&\te&&\Tor^R_i(M,N)=0\text{\it\ for all $i\gg 0$ implies\ }
\Ext_R^i(M,N)=0\text{\it\ for all $i\gg 0$},\\
&\et&&\Ext_R^i(M,N)=0\text{\it\ for all $i\gg 0$ implies\ }
\Tor^R_i(M,N)=0\text{\it\ for all $i\gg 0$},\\
\end{align*}
where $M$ and $N$ range over all finitely generated $R$-modules.
In addition, $\gap$ is the property that if $n$ consecutive
$\Ext_R^i(M,N)$ vanish, then $\Ext_R^i(M,N)$ vanishes
for all $i>\dim R$, where $n$ is a fixed positive 
integer depending only on $R$. 

Let $\gor$, respectively $\ci$, denote the property that the ring is 
Gorenstein, respectively a complete intersection. The known implications 
among these properties are summarized
by a diagram in \cite{JS1}, which we reproduce here.  There is
one added improvement, which represents one
of the main results of this paper, namely the irreversibility of the
rightmost implication (i.e., Corollary \ref{eecor}).


$$
\xymatrixrowsep{1pc}
\xymatrixcolsep{.31pc}
\xymatrix{
& &\qquad\et{\ \& \ }\te\ar@{<=>}[rrrr]
& & & &\te
\ar@{=>}[drrrr]_{(6)}& & & & & & & \\
\ci\ar@{<=}[drr]<-0.7ex>|-{\boldsymbol\vert}\ar@{=>}[urr]<-0.7ex>
\ar@{<=}[urr]<0.7ex>|-{\boldsymbol\vert} & & & & &
& & & & &\ee\ar@{<=}[rrrr]<0.7ex>|-{\boldsymbol\vert}
\ar@{=>}[rrrr]<-0.7ex>& & &
&\gor\\
& &\gap\ar@{<=}[ull]<-0.7ex>\ar@{=>}[rrrr]^{{(4)}} & & &
&
\ab\ar@{=>}[urrrr]^{(5)} & & & & & & &}
$$



We do not know whether the implications (4), (5) and (6) are
reversible, or whether there exist other implications. Thus, while we know that the classes $\gap$, $\ab$, $\te$ and $\ee$
are all strictly  contained  between complete intersections and Gorenstein
rings, we do not know whether they represent four {\it distinct} classes or
fewer. 
\end{remark}

\appendix
\section{}

Let $P$ and $R$ be as defined in \ref{define}. The purpose of this appendix 
is to prove Proposition \ref{prop}. We recall its statement below:

\begin{proposition} 
\label{prop1}
The ring $R$ is local, with maximal ideal $\fm$, 
and satisfies the following properties: 
\begin{enumerate}[\quad\rm(1)]
\item $R$ has Hilbert series $H_R(t)=1+6t+6t^2+t^3$. 
More precisely, a basis of $R$ over $k$ is given by the following 
fourteen elements:
       $$
1,\, t, \, u,\, v,\, x,\, y,\, z,\, tv,\,uv,\,vx,\,vy,\,vz,\,tx,\,vtx\,
$$
\item $R$ is Gorenstein, with $\Soc(R)=(tvx)$. 
\item $R$ is a Koszul algebra. 
\end{enumerate}
\end{proposition}

\begin{chunk}
Examining the defining ideal $I$ of $R$, we see that $\fm^2$ is a 
$6$-dimensional vector space with basis
$tv,uv,vx,vy,vz,tx$. Set $s=tvx$. We provide below a relevant part of the
multiplication table of $R$:

\begin{center}
\label{table}
\renewcommand{\arraystretch}{1.5}
\begin{tabular}{c| c c c c c c c c }
$\cdot$  &$tv$ &$uv$ &$vx$ &$vy$ &$vz$ &$tx$ \\
\hline
$t$   &0  &0   &$s$  &0  &0   &0  \\
$u$   &0 &0  &0  &0 &$s$ &0 \\
$v$   &0 &0  &0 &0  &0 &$s$\\
$x$   &$s$ &0 &$s$ &0 &0 &0 \\
$y$    &0 &0 &0 &$s$ &0 &0 \\
$z$   &0 &$s$ &0 &0 &0 &$-\alpha s$\\
\end{tabular}
\end{center}
It is clear from the table that $\fm^4=0$, and $\fm^3=sR$. If we can prove that
$s\ne 0$, then the table above shows that  the socle
of $R$ is $1$-dimensional, hence $R$ is Gorenstein. 
\end{chunk}

In order to prove that $s\ne 0$ we will use the
well-known fact that every Gorenstein quotient of $P$ with socle in
degree 3 corresponds uniquely to a form F of degree $3$ in a divided
power algebra \cite[2.12]{IK}. 
We recall below the necessary details of this fact.

\begin{chunk}
\label{divided}
Let $\mathcal D$ be the divided power algebra in $6$ variables 
{\small $\T_T, \T_U,\T_V,\T_X,\T_Y ,\T_Z$} over $k$.  
We refer to Appendix A in \cite{IK} for basic properties of this
algebra. If $M$ is the multi-index $M=(m_T,m_U,\dots, m_Z)\in \mathbb N^6$, 
then the divided powers monomials 
{\small $\T^{[M]}=\T_T^{[m_T]}\dots \T_Z^{[m_Z]}$} with $|M|=m_T+\cdots+m_Z=j$
form a $k$-basis of $\mathcal D_j$.  By definition  {\small $\T_T^{[1]}=\T_T$}, 
and similarly for the other variables.  Multiplication in this algebra
is defined by extending  by linearity the rule  
{\small
\[
\T^{[M]}\cdot\T^{[N]}=\frac{(M+N)!}{M!N!}\T^{[M+N]},
\]}
where {\small
${\displaystyle\frac{(M+N)!}{M!N!}=\frac{(m_T+n_T)!\cdots(m_Z+n_Z)!}{m_T!\cdots
m_Z!n_T!\cdots n_Z!}}$}. 

One defines an action of the polynomial ring $P$ on
$\mathcal D$ as follows:
If $M=(m_T,\dots m_Z)$ and $N=(n_T,\dots, n_Z)$ are non-negative
multi-indices in $\mathbb N^6$, then:
\[
(T^{m_T}U^{m_U}V^{m_V}X^{m_X}Y^{m_Y}Z^{m_Z})\circ \T^{[A]}=
\begin{cases}
\T^{[N-M]}\,\,\,  \text{if} \,\,M\le N\\
0\,\,\,\,\, \text{otherwise}
\end{cases}
\]
and then extend by linearity. 
 
Given a form $F\in\mathcal D$ of degree $j$, one can define the ideal
$I_F\subseteq P$ of the elements  in $P$ which annihilate $F$:
$$
I_F=\{G\in P\mid G\circ F=0\}
$$
It is well-known, cf. \cite[Section 2.3]{IK}, 
that the ring $P/I_F$ is a Gorenstein Artinian ring with socle in degree $j$. 
\end{chunk}

\begin{proof}[Proof of Proposition {\rm \ref{prop1}}]
Consider the ring $\mathcal D$ as in \ref{divided} and let
$F\in\mathcal D$ be the following form (or $1/6$ of this form if 
the characteristic of $k$ is 2 or 3):
{\small 
\[
-3\T_Z\T_Y^2 + \T_X^3-3\T_Z\T_T^2+3\T_Y\T_T^2 + 6\T_Z\T_U\T_V +  
3\T_Y^2\T_V + 3\T_X^2\T_V + 6\T_X\T_T\T_V-3\alpha\T_Z(\T_X+\T_T)^2.
\]}
It is easy to check that the generators of the ideal $I$ annihilate $F$, hence $I\subseteq I_F$. By \ref{divided},
the ring $P/I_F$ is Gorenstein and has socle in degree $3$. Consequently, 
we have  $s\ne 0$, hence $R$ is
Gorenstein with socle $sR$.  This proves (1) and (2). 

(3) In \ref{define} we have written the generators of $I$  
so that the first monomial occurring
in each generator is its initial term with respect to reverse
lexicographic order associated to the 
variable ordering $Z>U>Y>X>T>V$. Let $J$ denote the
ideal generated by these initial terms:
$$
J=(Z^2,\,UZ,\,U^2,\,YZ,\,UY,\,
Y^2,
XZ,\,UX,\,XY,\,X^2,\,
TZ,\,TU,
TY,\,T^2,\,V^2)\,.
$$

It is easy to  check that the Hilbert
series of $P/J$ is equal to $1+6t+6t^2+t^3$, and hence it
is equal to the Hilbert series of $R$.  It follows that the initial ideal of $I$ equals $J$, 
and this shows that the generators of $I$ listed above are a 
Gr\"obner basis for $I$. By \cite{F} this shows that $R$ is a Koszul
algebra. 
\end{proof}

\section*{Acknowledgments} The authors wish to thank Luchezar Avramov
for his valuable comments and suggestions, Kristen Beck and Srikanth 
Iyengar for work and suggestions on an earlier version of this paper.
They would also like to thank the referees for their helpful suggestions.

\end{document}